# TRANSCIENCE/RECURRENCE FOR NORMALLY REFLECTED BROWNIAN MOTION IN UNBOUNDED DOMAINS

By Ross G. Pinsky

*Technion—Israel Institute of Technology*

Let $D \subsetneq R^d$ be an unbounded domain and let $B(t)$ be a Brownian motion in $D$ with normal reflection at the boundary. We study the transcience/recurrence dichotomy, focusing mainly on domains of the form $D = \{(x,z) \in R^{l+m} : |z| < H(|x|)\}$, where $d = l + m$ and $H$ is a sufficiently regular function. This class of domains includes various horn-shaped domains and generalized slab domains.

**1. Introduction and statement of results.** Let $D \subsetneq R^d$ be an unbounded domain with $C^2$-boundary and consider Brownian motion $\{B(t), 0 \le t < \infty\}$ in $D$ with normal reflection at the boundary [4]. The normally reflected Brownian motion is a reversible Markov process with self-adjoint generator $\frac{1}{2}\Delta_N$, where $\Delta_N$ is the Neumann Laplacian on $D$. In this paper we study the transience/recurrence dichotomy for this process.

We recall the following facts which can be found in [3]. For a bounded open set $U \subset D$, let $\tau_U = \inf\{t \ge 0 : B(t) \in \bar{U}\}$ denote the first hitting time of $\bar{U}$, the closure of $U$. Let $P_y$ and $E_y$ denote probabilities and expectations for the normally reflected Brownian motion starting from $y \in \bar{D}$. One has either $P_y(\tau_U < \infty) = 1$, for all $y \notin \bar{U}$ and all bounded open $U$ for which $D - U$ is connected, or else, one has $P_y(\tau_U < \infty) < 1$, for all such $y$ and $U$. In the former case, the process is *recurrent* and in the latter case it is *transient* and one has $\lim_{t \to \infty} |B(t)| = \infty$ a.s. In the recurrent case, the process is called *positive recurrent* or *null recurrent*, respectively, according to whether it does or does not possess an invariant probability measure. One has either $E_y \tau_U < \infty$, for all $y$ and $U$ as above, or else one has $E_y \tau_U = \infty$, for all such $y$ and $U$. The former case is equivalent to positive recurrence.

The power of self-adjointness renders disarmingly simple the characterization of positive recurrence for the normally reflected Brownian motion:

PROPOSITION 1. *The normally reflected Brownian motion in $D \subset R^d$ is positive recurrent if and only if $\text{Vol}(D) < \infty$.*

PROOF. For one direction, note that, by self-adjointness, the transition density function $p(t, y, y')$ for the normally reflected Brownian motion is symmetric in $y$ and $y'$; therefore, $\int_D p(t, y, y') \, dy = \int_D p(t, y', y) \, dy = 1$. Thus, if $\text{Vol}(D) < \infty$, then the constant $\frac{1}{\text{Vol}(D)}$ is the invariant probability density. On the other hand, if there exists an invariant probability density $\phi$, then by the spectral theorem and the regularity, $\phi$ must be an eigenfunction corresponding to the eigenvalue 0 for $-\frac{1}{2}\Delta_N$; that is, $\phi$ is $\frac{1}{2}\Delta_N$-harmonic: $\frac{1}{2}\Delta_N \phi = 0$ in $D$ and $\nabla \phi \cdot n = 0$ on $\partial D$. The existence of the invariant probability measure guarantees recurrence, and it is well known that the generator of a recurrent Markov process possesses no nonconstant positive harmonic functions [3], Chapter 4; thus, we conclude that $\phi$ is constant. But if $\phi$ is constant and is a probability density on $D$, then $\text{Vol}(D) < \infty$. □

With regard to the transience/recurrence dichtomy, the self-adjointness no longer acts as a *deus ex machina*. However, it does allow one to prove that this dichotomy exhibits monotonicity with respect to the domain; that is, if the normally reflected Brownain motion in $D$ is recurrent, then the normally reflected Brownian motion in $\hat{D}$ is also recurrent if $\hat{D} \subset D$. The proof hinges on a variational formula which we now describe. Let $D_n = D \cap \{|y| < n\}$, assume without loss of generality that $D_1 = \{|y| < 1\}$ and let $A_{D,n}^{1,0} = \{u \in C^2(D_n - \bar{D}_1) \cap C(\bar{D}_n - D_1) : u = 1 \text{ on } \partial D_1 \text{ and } u = 0 \text{ on } \{|y| = n\} \cap D\}$. Then [2, 3]

$$(1.1) \quad l \equiv \lim_{n \to \infty} \inf_{u \in A_{D,n}^{1,0}} \int_{D_n} |\nabla u|^2 \, dy \begin{cases} = 0, & \text{if } B(t) \text{ is recurrent}; \\ > 0, & \text{if } B(t) \text{ is transient}. \end{cases}$$

If $\hat{D} \subset D$, then it is easy to see that any $u \in A_{D,n}^{1,0}$, when restricted to the appropriate smaller domain, belongs to $A_{\hat{D},n}^{1,0}$. The domain monotonicity follows from this and (1.1).

In particular, the domain monotonicity renders trivial the dichotomy of transience/recurrence when $d = 2$: since standard Brownian motion in all of $R^2$ is recurrent, normally reflected Brownian motion in any domain $D \subset R^2$ is a fortiori recurrent.

For the rest of this paper we will study the transience/recurrence dichotomy in domains $D \subset R^d$, $d \geq 3$, of the following form. Let $H > 0$ be a continuous function on $[0,\infty)$. For $l, m \geq 1$ with $l + m \geq 3$, denote points $y \in R^{l+m}$ by $y = (x, z)$, where $x \in R^l$ and $z \in R^m$. Let

(1.2) $$D = \{(x, z) \in R^{l+m} : |z| < H(|x|)\}.$$

We single out two particular subclasses of domains as above. In the case that $l = 1$, we will write $d - 1$ instead of $m$, and such domains will be called $d$-dimensional *horn-shaped domains*. In the case that $m = 1$, we will write $d - 1$ instead of $l$, and such domains will be called $d$-dimensional *generalized slab domains*.

Before presenting the results, we note the following point of departure. In consonance with the above notation, denote the normally reflecting Brownian motion as $B(t) = (X(t), Z(t))$. Consider the case that $H$ is constant. Then $D$ is the product of $R^l$ and a ball in $R^m$. (In particular, if $l = 1$, it is a cylinder and if $m = 1$, it is a slab.) Then $X(t)$ in isolation is a standard $l$-dimensional Brownian motion. If $l \leq 2$, then it is clear that the recurrence of $X(t)$ guarantees the recurrence of $B(t)$. By the domain monotonicity noted above, this recurrence will also hold for any domain that is contained in $D$. Thus, for $l \leq 2$, the interesting case occurs when $H$ is unbounded. On the other hand, if $l \geq 3$, then $X(t)$ is transient and, thus, so is $B(t)$. By domain monotonicity, this transience will also hold for any domain containing $D$. Thus, the interesting case occurs when $\lim_{s \to \infty} H(s) = 0$.

We will impose the following regularity conditions on $H$.

ASSUMPTION H. $H \in C^3([0, \infty))$ and satisfies the following conditions: $\lim_{s \to \infty} H(s)H''(s) = \lim_{s \to \infty} H'(s) = 0$, $H(s)H'(s) = o(s)$ as $s \to \infty$, $\int^\infty \frac{(H'(s))^3}{H(s)} ds < \infty$, $\int^\infty \frac{(H'(s))^2}{s} ds < \infty$.

REMARK. Assumption H allows for $H(s) \sim s^\gamma$ with $\gamma < 1$ but not with $\gamma \geq 1$. However, this growth restriction is in fact irrelevant, as will be noted in the remark following the theorem and examples below.

THEOREM 1. *Let $d = l + m \geq 3$ and let $D \subset R^d$ be the domain defined by (1.2) via a function $H$ which satisfies Assumption H.*

(i) *If $\int^\infty s^{1-l} H^{-m}(s) ds = \infty$ and $\hat{D} \subseteq D$, then normally reflected Brownian motion in $\hat{D}$ is recurrent.*

(ii) *If $\int^\infty s^{1-l} H^{-m}(s) ds < \infty$ and $D \subseteq \hat{D}$, then normally reflected Brownian motion in $\hat{D}$ is transient.*

REMARK. For a related result in the context of random walks, see [1], page 400.

EXAMPLE 1 (*Horn-shaped domains*). Denote points in $R^d$, $d \geq 3$, by $(x,z)$ with $x \in R$ and $z \in R^{d-1}$. Let $D$ be the $d$-dimensional horn-shaped domain defined by $D = \{(x,z) \in R^d : |z| < (1+|x|)^\gamma\}$. Then normally reflected Brownian motion in $D$ is recurrent if $\gamma \leq \frac{1}{d-1}$ and transient if $\gamma > \frac{1}{d-1}$.

EXAMPLE 2 (*Generalized slab domains*). Denote points in $R^d$, $d \geq 3$, by $(x,z)$ with $x \in R^{d-1}$ and $z \in R$. If $D$ is the $d$-dimensional generalized slab defined by $D = \{(x,z) \in R^d : |z| < (1+|x|)^\gamma\}$, then the normally reflected Brownian motion in $D$ is recurrent if $\gamma \leq 3-d$ and transient if $\gamma > 3-d$. If $d = 3$ and $D$ is the generalized slab domain defined by $D = \{(x,z) \in R^d : |z| < \log^\gamma(2+|x|)\}$, then the normally reflected Brownian motion in $D$ is recurrent if $\gamma \leq 1$ and transient if $\gamma > 1$.

EXAMPLE 3. For the domain $D = \{(x,z) \in R^{l+m} : |z| < H(|x|)\}$, with $H(s) \sim s^\gamma$, the normally reflected Brownian motion is recurrent if $\gamma \leq \frac{2-l}{m}$ and transient if $\gamma > \frac{2-l}{m}$.

REMARK. Example 3 demonstrates that for any domain $D$ defined by (1.2) with $H(s) \sim s^\gamma$, the value of $\gamma$ that represents the cut-off between recurrence and transience is never greater than $\frac{1}{2}$. This fact, along with the domain monotonicity with regard to transience and recurrence, renders irrelevant the growth condition imposed by Assumption H and noted in the remark following that assumption.

We now give the probabilistic intuition for Theorem 1. Consider first an $(l+m)$-dimensional Brownian motion $(X(t), Z(t))$ in the cylinder $\{(x,z) \in R^{l+m} : |z| < a\}$ with *skew-reflection* at $|z| = a$ corresponding to the reflection vector $(\gamma \frac{x}{|x|}, -\beta \frac{z}{|z|})$, where $\beta > 0$ and $\gamma$ are constants and the vector has been normalized by the requirement that $\gamma^2 + \beta^2 = 1$. Let $\rho(t) = |X(t)|$. By Itô's formula for diffusions with reflection, one has $\rho(t) = W(t) + \frac{l-1}{2}\int_0^t \frac{1}{\rho(s)}\,ds + \gamma L_a(t)$, where $W(t)$ is a one-dimensional Brownian motion and $L_a(t)$ is the local time up to time $t$ of $Z(\cdot)$ at $\{|z| = a\}$ [4]. We will show below that $E_{x,z} L_a(t) \sim \frac{m}{2a} t$ as $t \to \infty$. Thus, over the long run, the local time term behaves like a constant drift of strength $\frac{m}{2a}$. Consequently, the behavior of $\rho(t) = |X(t)|$ over the long run should be like the behavior of the one-diffusion generated by $\frac{1}{2}\frac{d^2}{d\rho^2} + \frac{l-1}{2\rho}\frac{d}{d\rho} + \frac{m\gamma}{2a}\frac{d}{d\rho}$.

With the above analysis in mind, consider the normally reflected Brownian motion in $D$, where $D$ is given by (1.2). At a point $(x,z) \in \partial D$, one has $|z| = H(|x|)$ and the normalized normal reflection vector pointing into $D$ is given by $((H')^2(|x|) + 1)^{-1/2}(H'(|x|)\frac{x}{|x|}, -\frac{z}{|z|})$. Taking our cue from the case above with constant skew reflection vector, and letting $a = H(|x|)$ and $\gamma = \frac{H'(|x|)}{((H')^2(|x|)+1)^{1/2}}$, we would expect the behavior of $|X(t)|$ over the long run to be like that of the one-dimensional diffusion generated by $\frac{1}{2}\frac{d^2}{d\rho^2} + \frac{l-1}{2\rho}\frac{d}{d\rho} + \frac{mH'(\rho)}{2H(\rho)((H')^2(\rho)+1)^{1/2}}\frac{d}{d\rho}$. Since we are assuming that $H'(\rho) \to 0$ as $\rho \to \infty$, the above diffusion is not much different from the one generated by $A \equiv \frac{1}{2}\frac{d^2}{d\rho^2} + \frac{l-1}{2\rho}\frac{d}{d\rho} + \frac{mH'(\rho)}{2H(\rho)}\frac{d}{d\rho}$. The criterion in Theorem 1 is exactly the transience/recurrence criterion [3], Chapter 5, for the one-dimensional diffusion whose generator is $A$.

We now return to show that $E_{x,z}L_a(t) \sim \frac{m}{2a}t$ as $t \to \infty$. The process $Z(t)$ in isolation is a Brownian on the $a$-ball with normal reflection at the boundary. Let $\mathcal{E}_z$ denote the expectation for this process starting from $z$. Let $\Delta_m$ denote the Laplacian on $R^m$, let $r = |z|$ and note that $\frac{1}{2}\Delta_m(|z|^2) = m$ and $\frac{\partial(|z|^2)}{\partial r}|_{|z|=a} = 2a$. Let $\sigma_r = \inf\{t \geq 0 : |Z(t)| = r\}$ and let $z_r$ denote a point satisfying $|z_r| = r$. We now twice apply Itô's formula for diffusions with reflection to the process $Z(t)$ and the function $|z|^2$ and take expectations. One time we start at $z_{a/2}$ and terminate at $\sigma_a$, and the other time we start at $z_a$ and terminate at $\sigma_{a/2}$. We obtain

$$a^2 = \mathcal{E}_{z_{a/2}}|Z(\sigma_a)|^2 = \frac{a^2}{4} + \mathcal{E}_{z_{a/2}}\int_0^{\sigma_a} m\,dt = \frac{a^2}{4} + m\mathcal{E}_{z_{a/2}}\sigma_a$$

and

$$\frac{a^2}{4} = \mathcal{E}_{z_a}|Z(\sigma_{a/2})|^2 = a^2 + \mathcal{E}_{z_a}\int_0^{\sigma_{a/2}} m\,dt - 2a\mathcal{E}_{z_a}L_a(\sigma_{a/2})$$
$$= a^2 + m\mathcal{E}_{z_a}\sigma_{a/2} - 2a\mathcal{E}_{z_a}L_a(\sigma_{a/2}).$$

Of course, we also have $\mathcal{E}_{z_{a/2}}L_a(\sigma_a) = 0$ since $L_a(t) = 0$ for $t \leq \sigma_a$. From the above equations we obtain

$$\frac{\mathcal{E}_{z_{a/2}}L_a(\sigma_a) + \mathcal{E}_{z_a}L_a(\sigma_{a/2})}{\mathcal{E}_{z_{a/2}}\sigma_a + \mathcal{E}_{z_a}\sigma_{a/2}} = \frac{m}{2a}.$$

Thus, $\frac{m}{2a}$ is equal to the expected local time accrued by $Z(t)$ at $\{|z| = a\}$ over a complete cycle starting at $|z| = \frac{a}{2}$ and ending again at $|z| = \frac{a}{2}$ after reaching $|z| = a$ divided by the expected total time of such a complete cycle. Using this, it follows easily from the renewal theorem that $E_{x,z}L_a(t) \sim \frac{m}{2a}t$ as $t \to \infty$.

The rest of the paper is organized as follows. In Section 2 we construct a class of functions from which we will cull appropriate Lyapunov functions which will be used to prove Theorem 1. Then in Section 3 we give the proof of Theorem 1.

**2. Lyapunov functions.** In this section we construct a family of functions from which we will cull appropriate Lyapunov functions which will be used to prove Theorem 1. Let $D$ be a domain in $R^{l+m}$ defined by (1.2) for some choice of $H$ satisfying Assumption H. Recall that we are writing points in $R^{l+m}$ by $(x, z)$, where $x \in R^l$ and $z \in R^m$. Let $\rho = |x|$ and $r = |z|$. The functions we construct will be functions of $\rho$ and $r$ and will denoted by $u(\rho, r)$. Each such function $u$ will depend on a function $f$, which is an arbitrary $C^2$-function from $[0, \infty)$ to $(0, \infty)$. For such an $f$, we define $u$ by

$$(2.1) \quad u\left(-\frac{1}{2}\frac{H'(s)}{H(s)}r^2 + s + \frac{1}{2}H(s)H'(s), r\right) = f(s), \qquad 0 \le r \le H(s).$$

Assumption H guarantees that, for sufficiently large $s$, the first argument in $u$ above is increasing in $s$ for each fixed $r \in [0, H(s)]$. Thus, $u$ is well defined by (2.1) for all sufficiently large $s$ and $r \le H(s)$. In particular then, for some $\rho_0 > 0$, $u$ is well defined on $\bar{D} \cap \{|x| \ge \rho_0\}$. Note that for each fixed $s$, the argument of $u$ in (2.1) defines a paraboloid with independent variable $r = |z|$. These paraboloids are the level sets of $u$. Note that they have been chosen so that

$$(2.2) \qquad (\nabla u \cdot n)(x, z) = 0 \qquad \text{for } (x, z) \in \partial D \cap \{|x| \ge \rho_0\},$$

where $n$ is the unit inward normal at $\partial D$.

We now calculate $\frac{1}{2}\Delta u$. Since $u$ depends only on $\rho = |x|$ and $r = |z|$, we have

$$(2.3) \qquad \frac{1}{2}\Delta u = \frac{1}{2}u_{\rho\rho} + \frac{l-1}{2\rho}u_\rho + \frac{1}{2}u_{rr} + \frac{m-1}{2r}u_r.$$

In the calculations that follow, we simplify notation by defining

$$L(s) \equiv \log H(s) \quad \text{and} \quad Q(s) \equiv H^2(s).$$

Also, we suppress the argument $(-\frac{1}{2}\frac{H'(s)}{H(s)}r^2 + s + \frac{1}{2}H(s)H'(s), r)$ that appears in $u$ and its derivatives. Differentiating (2.1) with respect to $r$ gives

$$(2.4) \qquad -rL'(s)u_\rho + u_r = 0$$

and differentiating a second time with respect to $r$ gives

$$(2.5) \qquad r^2(L'(s))^2 u_{\rho\rho} - L'(s)u_\rho + u_{rr} - 2rL'(s)u_{\rho r} = 0.$$

Differentiating (2.1) with respect to $s$ gives

$$(2.6) \qquad (-\tfrac{1}{2}r^2 L''(s) + 1 + \tfrac{1}{4}Q''(s))u_\rho = f'(s)$$

and differentiating a second time with respect to $s$ gives

$$(2.7) \quad (-\tfrac{1}{2}r^2 L''(s) + 1 + \tfrac{1}{4}Q''(s))^2 u_{\rho\rho} + (-\tfrac{1}{2}r^2 L'''(s) + \tfrac{1}{4}Q'''(s))u_\rho = f''(s).$$

Differentiating (2.4) with respect to $s$ gives

$$(-\tfrac{1}{2}r^2 L''(s) + 1 + \tfrac{1}{4}Q''(s))u_{\rho r}$$
$$- rL'(s)(-\tfrac{1}{2}r^2 L''(s) + 1 + \tfrac{1}{4}Q''(s))u_{\rho\rho} - rL''(s)u_\rho = 0. \tag{2.8}$$

Equation (2.4) allows one to solve for $u_r$ in terms of $u_\rho$ and $L$. Equation (2.8) allows one to solve for $u_{\rho r}$ in terms of $u_{\rho\rho}, u_\rho, L$ and $Q$, and then (2.5) allows one to solve for $u_{rr}$ in terms of these same functions. Equation (2.6) allows one to solve for $u_\rho$ in terms of $f, L$ and $Q$, and then (2.7) allows one to solve for $u_{\rho\rho}$ in terms of these same functions. We are therefore able to express $\tfrac{1}{2}\Delta u$ as given in (2.3) in terms of $f, L$ and $Q$ as follows:

$$\tfrac{1}{2}\Delta u = \tfrac{1}{2}(\Delta u)\left(-\tfrac{1}{2}\frac{H'(s)}{H(s)}r^2 + s + \tfrac{1}{2}H(s)H'(s), r\right)$$
$$= \tfrac{1}{2}A(r,s)f''(s) + \tfrac{1}{2}B(r,s)f'(s), \tag{2.9}$$

where

$$A(r,s) = \frac{1 + r^2(L'(s))^2}{C^2} \tag{2.10}$$

and

$$B(r,s) = m\frac{L'(s)}{C} + (l-1)\frac{1}{C\int C} + \frac{2r^2 L' L''}{C^2}$$
$$- \frac{(1 + r^2(L'(s))^2)(\partial C/\partial s)}{C^3} \tag{2.11}$$

with

$$C = C(r,s) = 1 + \tfrac{1}{4}Q''(s) - \tfrac{1}{2}r^2 L''(s) \tag{2.12}$$

and

$$\int C = \left(\int C\right)(r,s) = s + \tfrac{1}{4}Q'(s) - \tfrac{1}{2}r^2 L'(s). \tag{2.13}$$

We note that, by the assumptions on $H$, $C(r,s)$ and $(\int C)((r,s)$ are bounded away from 0 for large $s$ and $r \leq H(s)$.

**3. Proof of Theorem 1.** It suffices to prove the results for $D$ as in the statement of the theorem. The corresponding results for $\hat{D}$ then follow from the domain monotonicity property noted in Section 1. Recall that we are writing the normally reflected Brownian motion as $B(t) = (X(t), Z(t))$. For $x_0 > 0$, let $\tau_{x_0} = \inf\{t \geq 0 \colon X(t) \leq x_0\}$. For a $C^2$-function $v$ defined on $\bar{D} \cap \{x \geq x_0\}$, Itô's formula gives

$$\begin{aligned}
&v(B(t \wedge \tau_{x_0})) \\
&= v(B(0)) + \int_0^{t \wedge \tau_{x_0}} (\nabla v)(B(s)) \, dW(s) \\
&\quad + \int_0^{t \wedge \tau_{x_0}} \frac{1}{2}(\Delta v)(B(s)) \, ds + \int_0^{t \wedge \tau_{x_0}} (\nabla v \cdot n)(B(s)) \, dL(s),
\end{aligned} \tag{3.1}$$

where $W$ is a standard $d$-dimensional Brownian motion, $n$ is the inward unit normal at $\partial D$ and $L$ is the local time of $B(t)$ at the boundary $\partial D$. If $\nabla v \cdot n = 0$ on $\partial D \cap \{x \geq x_0\}$, then the local time term vanishes. If one can find such a $v$ for which $\Delta v \leq 0$ on $\bar{D} \cap \{x \geq x_0\}$ and such that $\lim_{x \to \infty} v(x, z) = \infty$, then the normally reflected Brownian motion is recurrent, while if one can find such a bounded $v$ for which $\Delta v \leq 0$ on $\bar{D} \cap \{x \geq x_0\}$ and such that $v(x_1, z) < \inf_{|z| \leq H(x_0)} v(x_0, z)$, for some $x_1 > x_0$, then the normally reflected Brownian motion is transient [3], Chapter 6.

Recall from (2.2) that the function $u$ defined in (2.1) satisfies $\nabla u \cdot n = 0$ on $\partial D \cap \{|x| \geq \rho_0\}$. If one can choose the function $f$ appearing in the definition of $u$ such that $\lim_{s \to \infty} f(s) = \infty$ and $\frac{1}{2}\Delta u(x, z) \leq 0$, for $(x, z) \in D$ with $x$ sufficiently large, then one can choose $v = u$ and conclude that the normally reflected Brownian motion is recurrent. Alternatively, if one can choose a bounded, increasing $f$ such that $\frac{1}{2}\Delta u(x, z) \geq 0$ for $(x, z) \in D$, with $x$ sufficiently large, then one can choose $v = -u$ and conclude that the normally reflected Brownian motion is transient.

From (2.9) it follows that if $\Gamma^+(s), \Gamma^-(s)$ satisfy

$$\Gamma^+(s) \geq \sup_{r \leq H(s)} \frac{B}{A}(r, s), \qquad \Gamma^-(s) \leq \inf_{r \leq H(s)} \frac{B}{A}(r, s) \tag{3.2}$$

and if we let

$$f^{\pm}(s) = \int_{s_0}^s dt \exp\left(-\int_{s_0}^t \Gamma^{\pm}(\rho) \, d\rho\right) \tag{3.3}$$

for some $s_0 > 0$, then $\frac{1}{2}\Delta u \leq 0$ if $u$ is constructed from $f^+$, and $\frac{1}{2}\Delta u \geq 0$ if $u$ is constructed from $f^-$. Thus, if one can choose $\Gamma^+$ so that $\lim_{s \to \infty} f(s) = \infty$, then one will obtain recurrence, while if one can choose $\Gamma^-$ so that $f^-$ is bounded, then one will obtain transience. The assumptions on $H$ guarantee that, in the calculations that follow below, the denominators of all the fractions as well as the numerators of the fractions that appear as arguments of logarithms are bounded away from 0 for $s$ sufficiently large. We choose $s_0$ above sufficiently large to accommodate this condition.

In order to choose $\Gamma^\pm$ appropriately, we must analyze $\frac{B}{A}$ from (2.10) and (2.11). After performing some algebra to isolate terms that do not depend on $r$, and arranging certain other terms as logarithmic derivatives, we obtain

$$\begin{aligned}(3.4) \quad \frac{B}{A}(r,s) = & \, mL'(s) + (l-1)\frac{(s + 1/4Q'(s) - 1/2r^2 L'(s))'}{s + 1/4Q'(s) - 1/2r^2 L'(s)} \\ & + \frac{1}{4}mL'(s)Q''(s) \\ & - mL'(s)\frac{r^2(L'(s))^2(1 + 1/4Q''(s))}{1 + r^2(L'(s))^2} \\ & + \frac{5-m}{2}\frac{(1 + r^2(L'(s))^2)'}{1 + r^2(L'(s))^2} \\ & - \frac{(1 + 1/4Q''(s) - 1/2r^2 L''(s))'}{1 + 1/4Q''(s) - 1/2r^2 L''(s)} \\ & - \frac{l-1}{2}\frac{r^2(L'(s))^2(1 + 1/4Q''(s) - 1/2r^2 L''(s))}{s + 1/4Q'(s) - 1/2r^2 L'(s)},\end{aligned}$$

where the prime denotes differentiation with respect to $s$. For fixed $s$, as functions of $r \in [0, H(s)]$, the three functions $\frac{(1+r^2(L'(s))^2)'}{1+r^2(L'(s))^2} = \frac{2r^2 L'(s) L''(s)}{1+r^2(L'(s))^2}$, $-\frac{(1+1/4Q''(s) - 1/2r^2 L''(s))'}{1+1/4Q''(s) - 1/2r^2 L''(s)} = -\frac{1/4Q'''(s) - 1/2r^2 L'''(s)}{1+1/4Q''(s) - 1/2r^2 L''(s)}$ and $\frac{(s+1/4Q'(s) - 1/2r^2 L'(s))'}{s+1/4Q'(s) - 1/2r^2 L'(s)} = \frac{1+1/4Q''(s) - 1/2r^2 L''(s)}{s+1/4Q'(s) - 1/2r^2 L'(s)}$ take on their maximum values and their minimum values at the endpoints $\{0, H(s)\}$. In particular then, the functions

$$(3.5) \quad \begin{aligned} E^+(s) &= \sup_{0 \le r \le H(s)} \frac{(1 + r^2(L'(s))^2)'}{1 + r^2(L'(s))^2}, \\ E^-(s) &= \inf_{0 \le r \le H(s)} \frac{(1 + r^2(L'(s))^2)'}{1 + r^2(L'(s))^2}, \end{aligned}$$

$$(3.6) \quad \begin{aligned} F^+(s) &= \sup_{0 \le r \le H(s)} \left(-\frac{(1 + 1/4Q''(s) - 1/2r^2 L''(s))'}{1 + 1/4Q''(s) - 1/2r^2 L''(s)}\right), \\ F^-(s) &= \inf_{0 \le r \le H(s)} \left(-\frac{(1 + 1/4Q''(s) - 1/2r^2 L''(s))'}{1 + 1/4Q''(s) - 1/2r^2 L''(s)}\right) \end{aligned}$$

and

$$(3.7) \quad \begin{aligned} G^+(s) &= \sup_{0 \le r \le H(s)} \frac{(s + 1/4Q'(s) - 1/2r^2 L'(s))'}{s + 1/4Q'(s) - 1/2r^2 L'(s)}, \\ G^-(s) &= \inf_{0 \le r \le H(s)} \frac{(s + 1/4Q'(s) - 1/2r^2 L'(s))'}{s + 1/4Q'(s) - 1/2r^2 L'(s)} \end{aligned}$$

are continuous and piecewise continuously differentiable, and

$$
\begin{aligned}
\int_{s_0}^t E^{\pm}(\rho)\,d\rho &= \log\frac{1+H_{E^{\pm}}^2(t)(L'(t))^2}{1+H_{E^{\pm}}^2(s_0)(L'(s_0))^2},\\
\int_{s_0}^t F^{\pm}(\rho)\,d\rho &= \log\frac{1+1/4Q''(s_0)-1/2H_{F^{\pm}}^2(s_0)(L''(s_0))^2}{1+1/4Q''(t)-1/2H_{F^{\pm}}^2(t)(L''(t))^2},\\
\int_{s_0}^t G^{\pm}(\rho)\,d\rho &= \log\frac{t+1/4Q'(t)-1/2H_{G^{\pm}}^2(t)L'(t)}{s_0+1/4Q'(s_0)-1/2H_{G^{\pm}}^2(s_0)L'(s_0)},
\end{aligned}
\tag{3.8}
$$

where $H_{*\pm}(s) \in \{0, H(s)\}$, for $* = E, F, G$.

Since $L' = \frac{H'}{H}$, $Q'' = 2HH'' + 2(H')^2$ and $0 \le r \le H$, it follows from Assumption H that, for some $C_0 > 0$, one has

$$
m\left|L'(s)\frac{r^2(L'(s))^2(1+1/4Q''(s))}{1+r^2(L'(s))^2}\right| \le C_0\frac{|H'(s)|^3}{H(s)}.
\tag{3.9}
$$

Similarly, for some $C_1 > 0$, one has

$$
\frac{l-1}{2}\frac{r^2(L'(s))^2|1+1/4Q''(s)-1/2r^2L''(s)|}{s+1/4Q'(s)-1/2r^2L'(s)} \le C_1\frac{(H'(s))^2}{s}.
\tag{3.10}
$$

In light of the above analysis, we define

$$
\begin{aligned}
\Gamma^+(s) &= mL'(s) + \frac{1}{4}mL'(s)Q''(s) + C_0\frac{|H'(s)|^3}{H(s)} + C_1\frac{(H'(s))^2}{s}\\
&\quad + \frac{5-m}{2}E^+(s) + F^+(s) + (l-1)G^+(s);\\
\Gamma^-(s) &= mL'(s) + \frac{1}{4}mL'(s)Q''(s) - C_0\frac{|H'(s)|^3}{H(s)} - C_1\frac{(H'(s))^2}{s}\\
&\quad + \frac{5-m}{2}E^-(s) + F^-(s) + (l-1)G^-(s).
\end{aligned}
\tag{3.11}
$$

These choices of $\Gamma^{\pm}$ satisfy (3.2). We define $f^{\pm}$ as in (3.3).

It remains to analyze the behavior of $f^{\pm}$ as $s \to \infty$. By (3.8) and the assumptions on $H$, $\exp(-\int_{s_0}^t (\frac{5-m}{2}E^{\pm}(\rho)+F^{\pm}(\rho))\,d\rho)$ is bounded and bounded from 0 for large $t$. By the assumptions on $H$, $\exp(\pm\int_{s_0}^t C_0\frac{|H'(\rho)|^3}{H(\rho)}\,d\rho)$ and $\exp(\pm\int_{s_0}^t C_1\frac{(H'(\rho))^2}{\rho}\,d\rho)$ are bounded and bounded from 0 for large $t$. We have

$$
\begin{aligned}
\int_{s_0}^t L'(\rho)Q''(\rho)\,d\rho &= 2\int_{s_0}^t \left(H'(\rho)H''(\rho) + \frac{(H'(\rho))^3}{H(\rho)}\right)d\rho\\
&= (H'(t))^2 - (H'(s_0))^2 + \int_{s_0}^t \frac{(H')^3}{H}(\rho)\,d\rho
\end{aligned}
$$

and, thus, by the assumptions on $H$, $\exp(-\frac{m}{4}\int_{s_0}^t L'(\rho)Q''(\rho)\,d\rho)$ is bounded and bounded from 0 for large $t$.

The assumptions on $H$ also show that

$$(1-c(s_0))\frac{t}{s_0} \leq \frac{t + 1/4Q'(t) - 1/2H_{G^\pm}^2(t)L'(t)}{s_0 + 1/4Q'(s_0) - 1/2H_{G^\pm}^2(s_0)L'(s_0)} \leq (1+c(s_0))\frac{t}{s_0}, \tag{3.12}$$

where $0 < c(s_0) = o(1)$ as $s_0 \to \infty$.

Thus, from (3.8) one has

$$\left(\frac{(1+c(s_0))t}{s_0}\right)^{1-l} \leq \exp\left(-\int_{s_0}^t (l-1)G^\pm(\rho)\,d\rho\right) \leq \left(\frac{(1-c(s_0))t}{s_0}\right)^{1-l}, \tag{3.13}$$

where $0 < c(s_0) = o(1)$ as $s_0 \to \infty$.

Finally, we have

$$\exp\left(-\int_{s_0}^t mL'(\rho)\,d\rho\right) = \left(\frac{H(t)}{H(s_0)}\right)^{-m}.$$

The calculations of the last two paragraphs along with (3.11) give

$$(1-C(s_0))\left(\frac{t}{s_0}\right)^{1-l}\left(\frac{H(t)}{H(s_0)}\right)^{-m} \leq \exp\left(-\int_{s_0}^t \Gamma^\pm(\rho)\,d\rho\right)$$

$$\leq (1+C(s_0))\left(\frac{t}{s_0}\right)^{1-l}\left(\frac{H(t)}{H(s_0)}\right)^{-m}, \tag{3.14}$$

where $0 < C(s_0) = o(1)$ as $s_0 \to \infty$.

We therefore conclude from (3.3) and (3.14) that $\lim_{s\to\infty} f^+(s) = \infty$, if $\int^\infty t^{1-l}H^{-m}(t)\,dt = \infty$, and $f^-$ is bounded, if $\int^\infty t^{1-l}H^{-m}(t)\,dt < \infty$.

DEPARTMENT OF MATHEMATICS
TECHNION—ISRAEL INSTITUTE OF TECHNOLOGY
HAIFA, 32000
ISRAEL
E-MAIL: pinsky@math.technion.ac.il
URL: http://www.math.technion.ac.il/~pinsky/